\documentclass [a4 paper,12pt]{article}
\usepackage{mathrsfs}

\usepackage[usenames]{color}
\usepackage{bbm,amsmath,amssymb,amsfonts,amsthm,graphics,graphicx}

\newtheorem{lem}{Lemma}
\newtheorem{thm}[lem]{Theorem}
\newtheorem{cor}[lem]{Corollary}
\newtheorem{con}{Conjecture}
\newtheorem{rem}{Remark}

\title{The asymptotic value of Randi\'{c} index for trees\footnote{Supported by
NSFC No.10831001, PCSIRT and the ``973" program.}}
\author{\small Xueliang Li, Yiyang Li\\
\small Center for Combinatorics and LPMC-TJKLC\\
\small Nankai University, Tianjin 300071, China}
\date{ }

\begin{document}

\maketitle

\begin{abstract}
Let $\mathcal{T}_n$ denote the set of all unrooted and unlabeled
trees with $n$ vertices, and $(i,j)$ a double-star. By assuming that
every tree of $\mathcal{T}_n$ is equally likely, we show that the
limiting distribution of the number of occurrences of the
double-star $(i,j)$ in $\mathcal{T}_n$ is normal. Based on this
result, we obtain the asymptotic value of Randi\'c index for trees.
Fajtlowicz conjectured that for any connected graph the Randi\'c
index is at least the average distance. Using this asymptotic value,
we show that this conjecture is true not only for almost all
connected graphs but also for almost all trees.\\[3mm]
\end{abstract}

\section{Introduction}

Let $\mathcal{T}_n$ denote the set of all unrooted and unlabeled
trees $T_n$ with $n$ vertices. A {\it pattern} $\mathcal{M}$ is a
given subtree. We say that $\mathcal{M}$ {\it occurs} in a tree if
$\mathcal{M}$ is a subtree of $T_n$ such that except for the
vertices of $\mathcal{M}$ with degree $1$, the other vertices must
have the same degrees with the corresponding vertices in $T_n$.
Surely, we can also let the vertices with degree $1$ match with each
other. Set $t_n=|\mathcal{T}_n|$. We introduce two functions:
$$t(x)=\sum_{n\geq1}t_nx^n,$$
$$t(x,u)=\sum_{n\geq1,k\geq 0}t_{n,k}x^nu^k,$$
where the coefficients $t_{n,k}$ denote the number of trees in
$\mathcal{T}_n$ that have $k$ occurrences of the pattern
$\mathcal{M}$. We assume that every tree of $\mathcal{T}_n$ is
equally likely. Let $X_n$ denote the number of occurrences of
$\mathcal{M}$ in a tree of $\mathcal{T}_n$. Therefore, $X_n$ is a
random variable on $\mathcal{T}_n$ with probability
$$\mbox{Pr}[X_n=k]=\frac{t_{n,k}}{t_n}.$$

In \cite{gk}, Kok showed that for any pattern $\mathcal{M}$ the
limiting distribution of $(X_n-EX_n)/$ $\sqrt{Var X_n}$ is a
distribution with density of the form $(A+Bt^2)\mbox{ exp}^{-Ct^2}$,
and $E(X_n)=(\mu+o(1)) n$ and $Var (X_n)=(\sigma+o(1))n$, where $A,
B, C, \mu, \sigma$ are some constants. Clearly, if $B=0$, it is a
normal distribution. It has been showed that if the pattern is a
{\it star} or a {\it path}, the corresponding distribution is
asymptotically normal. We refer the readers to \cite{gk, rs, dg} for
more details.

Recall that a {\it path} is a graph with a sequence of vertices with
an edge between every two consecutive vertices. A {\it star} is a
complete bipartite graph such that one partition contains only one
vertex, and we call this vertex the {\it center} of the star. A {\it
double-star} is a graph which is formed from two stars by connecting
their centers with an edge.

In this paper, we will show that if the pattern is a {\it
double-star}, the corresponding limiting distribution is also a
normal distribution, and get an estimate for the number $X_n$ of
occurrences of a double-star for almost all trees. Based on the
result, we then obtain the asymptotic value of {\it Randi\'{c}
index} for almost all trees in $\mathcal{T}_n$. The Randi\'c index
was introduced by Randi\'{c} \cite{r} in 1975, and later,
Bollob{\'{a}s and Erd\"{o}s \cite{be} generalized it to the {\it
general Randi\'{c} index}. The definition will be given in Section
3, and for a detailed survey we refer the readers to \cite{lishi}.

There is a conjecture on the relation between the Randi\'{c} index
and the average distance of a connected graph, proposed by
Fajtlowicz in \cite{fa}.
\begin{con}\label{rd}
Let $R(G)$ and $D(G)$ denote, respectively, the Randi\'c index and
the average distance of a graph $G$. Then, for any connected graph
$G$, $R(G)\geq D(G)$.
\end{con}
We will show that the conjecture is true not only for almost all
connected graphs but also for almost all trees.

In Section 2, we explore the limiting distribution of $X_n$
corresponding to a double-star. In Section 3, we apply the results
in Section 2 to considering the Randi\'{c} index.

\section{The distribution of $X_n$ for a double-star}

In this section, we concentrate on the limiting distribution of
$X_n$ for a double-star. Throughout this paper, we use $(i,j)$ to
denote the double-star with one vertex corresponding to a center of
degree $i$ and the other of degree $j$. Evidently, the number of
occurrence of $(i,j)$ in a tree is the number of edges in the tree
such that one end of the edge is of degree $i$ while the other is of
degree $j$. Without loss of generality, we always assume $i\leq j$.

In what follows, we first introduce some terminology and and
notation, which will be used in the sequel. For those not defined
here, we refer the readers to the book \cite{hp}.

Analogous to trees, we have the generating functions for rooted
trees and planted trees. Let $\mathcal{R}_n$ be the set of all
rooted trees with $n$ vertices and $r_n=|\mathcal{R}_n|$. We have
$$r(x)=\sum_{n\geq1}r_nx^n,$$
$$r(x,u)=\sum_{n\geq1,k\geq0}r_{n,k}x^nu^k,$$
and $r_{n,k}$ is the number of all rooted trees in $\mathcal{T}_n$
that have $k$ occurrences of $(i,j)$. A {\it planted tree} is
formed from a rooted tree and a new vertex by connecting the
vertex and the root of the rooted tree with a new edge. The new
vertex is called the {\it plant}, and we never count it in the
sequel. Let $\mathcal {P}_n$ denote the set of all planted trees
with $n$ vertices and $p_n=|\mathcal{P}_n|$. Then, we have
generating functions:
$$p(x)=\sum_{n\geq1}p_nx^n,$$
$$p(x,u)=\sum_{n\geq1,k\geq0}p_{n,k}x^nu^k,$$
where $p_{n,k}$ denote the number of planted trees in
$\mathcal{P}_n$ that have $k$ occurrences of $(i,j)$. By the
definitions of planted trees and rooted trees, it is easy to see
that
$$r(x,1)=r(x)=p(x,1)=p(x).$$
Furthermore, suppose the radius of convergence of $r(x)$ is $x_0$,
Otter \cite{ot} showed that $x_0$ satisfies $r(x_0)=1$ and the
asymptotic expansion of $r(x)$ is
\begin{eqnarray}\label{rexpa}r(x)=1-b(x_0-x)^{1/2}+c(x_0-x)
+d(x_0-x)^{3/2}+\cdots,\end{eqnarray}
where $x_0\approx 0.3383219$ and $b\approx 2.68112266$.

Let $\textbf{y}(x,u)=(y_1(x,u),\ldots,y_N(x,u))^T$ be a column
vector. We suppose $G(x,\textbf{y},u)$ is an analytic function with
non-negative Taylor coefficients. $G(x,\textbf{y},u)$ can be
expanded as
$$ G(x,\textbf{y},u)=\sum_{n,k}g_{n,k}x^nu^k.$$
Let $X_n$ denote a random variable with probability
\begin{equation}\label{rv}
\mbox{Pr}[X_n=k]=\frac{g_{n,k}}{g_n},
\end{equation}
where $g_n=\sum_kg_{n,k}$.

To show the limiting distribution of the number of occurrences of
the double-star $(i,j)$ for all trees is normal, we need a useful
lemma, which was used to explore the distribution of the number of
occurrences of a pattern for some other families of trees, such as
planar trees, labelled trees, rooted trees, {\it et al}. We refer
the readers to \cite{cdkk} for more details.

\begin{lem}\label{mainlem}
Let
$\textbf{F}(x,\textbf{y},u)=(F_1(x,\textbf{y},u),\ldots,F_N(x,\textbf{y},u))^T$
be functions analytic around $x=0$,
$\textbf{y}=(y_1,\ldots,y_N)^T=\textbf{0}$, $u=0$, with Taylor
coefficients that are all non-negative. Suppose
$\textbf{F}(0,\textbf{y},u)=\textbf{0}$,
$\textbf{F}(x,\textbf{0},u)\neq\textbf{0}$,
$\textbf{F}_x(x,\textbf{y},u)\neq \textbf{0}$, and for some $j$,
$\textbf{F}_{y_jy_j}(x,\textbf{y},u)\neq\textbf{0}$. Furthermore,
assume that $x=x_0$, $\textbf{y}=\textbf{y}_0$ is a non-negative
solution of the system of equations
\begin{align}
&\textbf{y}=\textbf{F}(x,\textbf{y},1)\\
&0=\mbox{det}(\textbf{I}-\textbf{F}_\textbf{y}(x,\textbf{y},1))\label{det}
\end{align}
inside the region of convergence of $\textbf{F}$, and $\textbf{I}$
is the unit matrix. Let $\textbf{y}=(y_1(x,u),\ldots,y_N(x,u))^T$
denote the analytic solution of the function system
\begin{equation}\label{fs}
\textbf{y}=\textbf{F}(x,\textbf{y},u)
\end{equation}
with $\textbf{y}(0,u)=\textbf{0}$. Moreover, let
$G(x,\textbf{y},u)$ be an analytic function with non-negative
Taylor coefficients such that the point
$(x_0,\textbf{y}(x_0,1),1)$ is contained in the region of
convergence. Finally, let $X_n$ be the random variable defined in
(\ref{rv}).

If the dependency graph $G_{\textbf{F}}$ of the function system
(\ref{fs}) is strongly connected, then the random variable $X_n$ is
asymptotically normal with mean
$$E(X_n)=\mu n+O(1)\mbox{  } (n\rightarrow \infty),$$
and variance
$$Var(X_n)=\sigma n+O(1)\mbox{ } (n\rightarrow \infty).$$

Moreover, suppose $\textbf{v}^T$ is the vector satisfying
$\textbf{v}^T(\textbf{I}-\textbf{F}_\textbf{y}(x_0,\textbf{y}_0,1))=0$,
we have that
\begin{eqnarray}\label{mu}
\mu=\frac{1}{x_0}\frac{\textbf{v}^T\textbf{F}_u(x_0,
\textbf{y}_0,1)}{\textbf{v}^T\textbf{F}_x(x_0, \textbf{y}_0,1)},
\end{eqnarray}
where $\textbf{F}_u$ and $\textbf{F}_x$ are the partial derivatives
of $\textbf{F}(x, \textbf{y},u)$.
\end{lem}

\begin{rem}
The {\it dependency graph} $G_{\textbf{F}}$ of
$\textbf{y}=\textbf{F}(x,\textbf{y},u)$ is strongly connected, if
there is no subsystem of equations that can be solved independently
from others. If $G_\textbf{F}$ is strongly connected, then
$\textbf{I}-\textbf{F}_\textbf{y}(x_0,\textbf{y}_0,1)$ has rank
$N-1$, i.e., $\textbf{v}$ is unique up to a nonzero factor. We refer
the readers to \cite{cdkk, d} for more details.
\end{rem}

Now, we consider the asymptotic distribution of $X_n$ corresponding
to the pattern $(i,j)$. Our main contribution is to establish
functional equations and apply Lemma \ref{mainlem} to obtain Theorem
\ref{mainthm}. For different $i,j$, we distinguish the following
three cases. Since only the tree with exactly two vertices contains
the pattern $(1,1)$, we do not need to consider the case for
$i=j=1$.

{\bf Case 1.} $i\neq j>1$.

We split $\mathcal{P}_n$ into three subsets according to the
degree of the root: the root is of degree $i$, $j$ and neither $i$
nor $j$, and we respectively let $a_i(x,u)$, $a_j(x,u)$ and
$a_0(x,u)$ be the generating functions (or $a_i$, $a_j$, $ a_0$
for short). It is easy to see that
\begin{eqnarray}\label{0ijp}
a_0(x,u)+a_i(x,u)+a_j(x,u)=p(x,u).
\end{eqnarray}

In what follows, there appears an expression of the form $Z(S_n,
f(x,u))$ (or $f(x)$), which is the substitution of the counting
series $f(x,u)$ (or $f(x)$) into the cycle index $Z(S_n)$ of the
symmetric group $S_n$. This involves replacing each variable $s_i$
in $Z(S_n)$ by $f(x^i,u^i)$ (or $f(x^i)$). For instance, if $n=3$,
then $Z(S_3)=(1/3!)(s_1^3+3s_1s_2+2s_3)$ and
$Z(S_3,f(x,u))=(1/3!)(f(x,u)^3+3f(x,u)f(x^2,u^2)+2f(x^3,u^3))$. We
refer the readers to \cite{hp} for details.

Employing the classical P\'{o}lya Enumeration Theorem, we have
$Z(S_{k-1};p(x))$ as the counting series of the planted trees whose
roots have degree $k$, and the coefficients of $x^p$ in
$Z(S_{k-1};p(x))$ is the number of planted trees of order $p+1$ (see
\cite{hp} p.51--54). Therefore, $p(x)$ satisfies
\begin{eqnarray}\label{exp}p(x)=x\sum_{k\geq0}Z(S_k;p(x))=x\mbox{
exp}^{\sum_{k\geq1}\frac{1}{k}p(x^k)}.
\end{eqnarray}
By the same way, we can obtain the following functional equations
\begin{align}
a_0(x,u)&=x \mbox{ exp}^{\sum_{k\geq1}\frac{1}{k}p(x^k,u^k)}-x
Z(S_{i-1};p(x,u))-x Z(S_{j-1};p(x,u)),\label{a0ij}\\
a_i(x,u)&=x{\sum_{\ell_1+\ell_2=i-1}Z(S_{\ell_1};
a_0(x,u)+a_i(x,u))\cdot
Z(S_{\ell_2}; a_j(x,u))}u^{\ell_2},\label{aiij}\\
a_j(x,u)&=x{\sum_{m_1+m_2=j-1}Z(S_{m_1}; a_0(x,u)+a_j(x,u))\cdot
Z(S_{m_2}; a_i(x,u))}u^{m_2}.\label{ajij}
\end{align}
For $a_0(x,u)$, since the degrees of the roots are neither $i$ nor
$j$, therefore there are two minor modifications. For $a_i(x,u)$,
if there exist $\ell_2$ vertices of degree $j$ adjacent to the
root, we should count $\ell_2$ occurrences of $(i,j)$ in addition,
and thus it is of $Z(S_{\ell_1}; a_0(x,u)+a_i(x,u))\cdot
Z(S_{\ell_2}; a_j(x,u))u^{\ell_2}$. Analogously, the equation of
$a_j(x,u)$ follows.

Then, for rooted trees, we have
\begin{eqnarray*}
r(x,u)&=&x\mbox{ exp}^{\sum_{k\geq1}\frac{1}{k}p(x^k,u^k)} -x Z(S_i;
p(x,u))-xZ(S_j; p(x,u))\\
&&+x\sum_{\ell_1+\ell_2=i}Z(S_{\ell_1}; a_0(x,u)+a_i(x,u))\cdot
Z(S_{\ell_2};
a_j(x,u))u^{\ell_2}\\
&&+x\sum_{m_1+m_2=j}Z(S_{m_1}; a_0(x,u)+a_i(x,u))\cdot Z(S_{m_2};
a_j(x,u))u^{m_2}
\end{eqnarray*}

In order to get the generating function for general trees, we need
the following lemma, which was used in \cite{ot} to get the famous
equation
\begin{eqnarray}\label{trp}
t(x)=r(x)-\frac{1}{2}p(x)^2+\frac{1}{2}p(x^2).
\end{eqnarray}
We can also obtain a similar equation for $t(x,u)$ from this
lemma.

Two edges in a tree are {\it similar}, if they are the same under
some automorphism of the tree. To {\it join} two planted trees is
to connect the two roots of the trees with a new edge and get rid
of the two plants. If the two planted trees are the same, we say
that the new edge is {\it symmetric}.
\begin{lem}\label{otter}
For any tree, the number of rooted trees corresponding to this tree
minus the number of nonsimilar edges (except the symmetric edge) is
the number $1$.
\end{lem}

Note that, if we delete any one edge of a similar set in a tree,
the yielded trees are the same two trees. Hence, different pairs
of planted trees correspond to nonsimilar edges. We refer the
readers to \cite{ot} for details. Then, analogous to (\ref{trp}),
we have
\begin{eqnarray}\label{tu}
t(x,u)=r(x,u)-\frac{1}{2}p(x,u)^2+\frac{1}{2}p(x^2,u^2)+a_i(x,u)\cdot
a_j(x,u)(1-u).
\end{eqnarray}
The last term serves to count the occurrences of $(i,j)$ when
joining two planted trees to form a tree, in which one has a root of
degree $i$ and the other has a root of degree $j$.

Now, we will use Lemma \ref{mainlem} to show that the distribution
of $X_n$ converges to a normal distribution and get the asymptotic
value of $E(X_n)$ corresponding to $(i,j)$.

We just need to verify that the system of functions (\ref{a0ij}),
(\ref{aiij}),(\ref{ajij}) satisfies equation (\ref{det}), since the
other conditions are easy to illustrated. We still denote the system
of functions by $\textbf{F}$. It is the function of vector
$\textbf{a}(x,u)=(a_0(x,u), a_i(x,u), a_j(x,u))^T$. Let
$\textbf{F}_{a_0}, \textbf{F}_{a_i}, \textbf{F}_{a_j}$ be the
partial derivations, respectively. Combining the fact that the
partial derivative enjoys (see \cite{dg})
$$\frac{\partial}{\partial s_1}Z(S_n;
s_1,\ldots,s_n)=Z(S_{n-1};s_1,\ldots,s_{n-1}),$$ with (\ref{rexpa}),
we obtain that
\begin{align*}\textbf{F}_{a_0}(x_0,\textbf{a}(x_0,1),1)&=\begin{pmatrix}
1-x_0Z(S_{i-2};p(x_0,1))-x_0Z(S_{j-2};p(x_0,1))
\\x_0\sum_{\ell_1+\ell_2=i-1}Z(S_{\ell_1-1};a_0(x_0,1)+a_i(x_0,1))Z(S_{\ell_2};a_j(x_0,1))
\\x_0\sum_{r_1+r_2=j-1}Z(S_{r_1-1};a_0(x_0,1)+a_j(x_0,1))Z(S_{r_2};a_i(x_0,1))
\end{pmatrix}\\
&=\begin{pmatrix} 1-x_0Z(S_{i-2};p(x_0,1))-x_0Z(S_{j-2};p(x_0,1))
\\x_0Z(S_{i-2};p(x_0,1))
\\x_0Z(S_{j-2};p(x_0,1))
\end{pmatrix}.\end{align*}
Similarly, we can get that
$\textbf{F}_{a_i}(x_0,\textbf{a}(x_0,1),1)=
\textbf{F}_{a_j}(x_0,\textbf{a}(x_0,1),1)=\textbf{F}_{a_0}(x_0,\textbf{a}(x_0,1),1).$
Therefore, one can readily see that
$$\mbox{det}(\textbf{I}-\textbf{F}_{\textbf{a}}(x_0,\textbf{a}(x_0,1),1))=0.$$
Moreover, from equation (\ref{trp}) it follows that
$$t(x_0,1)=(1+r(x_0^2))/2.$$
Note that $x_0<1$, and thus $x_0^2$ is surely inside the region of
convergence of $r(x)$. So, for the generating function $t(x,u)$, all
the conditions required by Lemma \ref{mainlem} are satisfied. Thus,
the distribution of $X_n$ corresponding to $(i,j)$, $i\neq j>1$, is
asymptotically normal.

From the form of $\textbf{F}_{\textbf{a}}(x_0,\textbf{a}(x_0,1),1)$,
it is not difficult to obtain that $\textbf{v}^T=(1,1,1)$ is a basic
solution. In what follows, we will compute
$\textbf{v}^T\textbf{F}_{x}(x_0,\textbf{a}(x_0,1),1)$ and
$\textbf{v}^T\textbf{F}_{u}(x_0,\textbf{a}(x_0,1),1)$ to estimate
$\mu$, which would be more brief than just to do with
$\textbf{F}_{x}(x_0,\textbf{a}(x_0,1),1)$ and
$\textbf{F}_{u}(x_0,\textbf{a}(x_0,1),1)$. Then, we have
\begin{align}
\textbf{v}^T\textbf{F}_{x}(x_0,\textbf{a}(x_0,1),1)&=\frac{1}{x_0}
+\sum_{k=2}p_x(x_0^k,1)x_0^{k-1},\label{fx}\\
\textbf{v}^T\textbf{F}_{u}(x_0,\textbf{a}(x_0,1),1)&=
\sum_{k=2}p_u(x_0^k,1)\nonumber\\
&+x_0\sum_{\ell_1+\ell_2=i-1}Z(S_{\ell_1};a_0(x_0,1)+a_i(x_0,1))Z(S_{\ell_2};a_j(x_0,1))\cdot
\ell_2\nonumber\\
&+x_0\sum_{m_1+m_2=j-1}Z(S_{r_1};a_0(x_0,1)+a_j(x_0,1))
Z(S_{m_2};a_i(x_0,1))\cdot m_2.\label{fu}
\end{align}
In view of $p(x,1)=p(x)=t(x)$, combining with (\ref{rexpa}) and
(\ref{exp}), it follows that
$$\frac{1}{x_0}+\sum_{k=2}p_x(x_0^k,1)x_0^{k-1}=\frac{p_x(x,1)(1-p(x,1))}{p(x,1)}{\Big
|}_{x=x_0}=b^2/2,$$
and thus
$$\textbf{v}^T\textbf{F}_{x}(x_0,\textbf{a}(x_0,1),1)=\frac{b^2}{2}.$$
However, we failed to do any further simplification for (\ref{fu}).
For convenience, denote the value of
$\textbf{v}^T\textbf{F}_{u}(x_0,\textbf{a}(x_0,1),1)$ by $w(i,j)$.
One can use a computer to get an approximate value of it. Thus,
$$\mu=\frac{2}{x_0b^2} w(i,j).$$

{\bf Case 2.} $i=1$, $j>1$.

We proceed to obtain the result in a same way as in Case 1. We
still use the same notation. But notice that when we split up
$\mathcal {P}_n$ according to the degrees of the roots, there
exists only one planted tree with root of degree $1$, that is, the
tree with only two nodes. Thus, we have
$$x+a_0(x,u)+a_j(x,u)=p(x,u),$$
and the system of functions is as follows
\begin{align}
a_0(x,u)&=x\mbox{
exp}^{\sum_{k\geq1}\frac{1}{k}p(x^k,u^k)}-x-xZ(S_{j-1};p(x,u)),\label{a01j}\\
a_j(x,u)&=x\sum_{m_1+m_2=j-1}Z(S_{m_1};p(x,u)-x)x^{m_2}u^{m_2}.\label{aj1j}
\end{align}

The same as previous, we can establish the generating function for
rooted trees
\begin{align*}r(x,u)&=x\mbox{
exp}^{\sum_{k\geq1}\frac{1}{k}p(x^k,u^k)}-xa_j(x,u)(1-u)\\
&-x\sum_{m_1+m_2=j} Z(S_{m_1};p(x,u)-x)\cdot
Z(S_{m_2};x)(1-u^{m_2}),
\end{align*}
and for general trees
\begin{align}t(x,u)&=r(x,u)-\frac{1}{2}p(x,u)^2+
\frac{1}{2}p(x^2,u^2)+xa_j(x,u)(1-u).\label{t1j}
\end{align}
It is not difficult to verify that (\ref{a01j}), (\ref{aj1j}) and
(\ref{t1j}) satisfy the conditions of Lemma \ref{mainlem}. We can
obtain $\textbf{v}^T=(1,1)$,
\begin{align*}
&\textbf{v}^T\textbf{F}_{x}(x_0,\textbf{a}(x_0,1),1)\\
&= x\mbox{ exp}^{\sum_{k\geq1}\frac{1}{k}p(x^k,u^k)} (1+\sum_{k\geq
2}p_x(x^k,u^k)x^{k-1})+\mbox{
exp}^{\sum_{k\geq1}\frac{1}{k}p(x^k,u^k)}-1\big|_{(x=x_0,u=1)}\\
&=\frac{b^2}{2},
\end{align*}
and
$$\textbf{v}^T\textbf{F}_{u}(x_0,\textbf{a}(x_0,1),1)=
\sum_{k\geq2}p_u(x_0^k,1)+x_0\sum_{\ell_1+\ell_2=j-1}Z(S_{\ell_1};
p(x_0,1)-x_0)x_0^{\ell_2}\cdot \ell_2.$$ Again, for convenience,
we denote $\textbf{v}^T\textbf{F}_{u}(x_0,\textbf{a}(x_0,1),1)$ by
$w(1,j)$. Then, it follows that $$\mu=\frac{2}{x_0b^2}w(1,j).$$

{\bf Case 3.} $i=j>1$.

Since the procedure is the same as previous, we leave out the
details of the proof for brevity. However, we still use the same
notation here without any conflicts.
\begin{align*}
a_0(x,u)&+a_j(x,u)=p(x,u),\\
 a_0(x,u)&=x\mbox{
exp}^{\sum_{k\geq1}\frac{1}{k}p(x^k,u^k)}-xZ(S_{j-1};p(x,u)),\\
a_j(x,u)&=x
\sum_{m_1+m_2=j-1}Z(S_{m_1};a_0(x,u))Z(S_{m_2};a_j(x,u))u^{m_2}.
\end{align*}

For general trees, we have
\begin{align*}t(x,u)&=r(x,u)-\frac{1}{2}p(x,u)^2+\frac{1}{2}p(x^2,u^2)\\
&+\frac{1}{2}a_j(x,u)\cdot
a_j(x,u)(1-u)-\frac{1}{2}a_j(x^2,u^2)(1-u).
\end{align*}

Employing Lemma \ref{mainlem}, asymptotic analysis of the functional
equations will give that $\textbf{v}^T=(1,1)$,
$\textbf{v}^T\textbf{F}_{x}(x_0,\textbf{a}(x_0,1),1)=b^2/2$ and
\begin{align*}
&\textbf{v}^T\textbf{F}_{u}(x_0,\textbf{a}(x_0,1),1)\\&
=\sum_{k\geq2}p_u(x_0^k,1)+x_0\sum_{m_1+m_2=j-1}Z(S_{m_1};a_0(x_0,1))Z(S_{m_2};a_j(x_0,1))\cdot
m_2.
\end{align*}

Then, we obtain that
$$\mu=\frac{2}{x_0b^2}w(j,j),$$
where $w(j,j)$ denote the value of
$\textbf{v}^T\textbf{F}_{u}(x_0,\textbf{a}(x_0,1),1).$

As a conclusion, we can establish the following theorem now.

\begin{thm}\label{mainthm}
Suppose $X_n$ is the random variable corresponding to the
occurrences of pattern $(i,j)$. The probability measure of $X_n$ is
defined as (\ref{rv}) for the generating function of trees $t(x,u)$.
Then, the distribution of $X_n$ is asymptotically normal with mean
$$EX_n=\frac{2}{x_0b^2}\cdot w(i,j)n+O(1)$$
and variance $Var X_n=\sigma(i,j)n+O(1)$, where $w(i,j)$,
$\sigma(i,j)$ are some constants.
\end{thm}

Following the book \cite{BB}, we will say that {\it almost every}
(a.e.) graph in a graph space $\mathcal{G}_n$ has a certain property
$Q$ if the probability $\mbox{Pr}(Q)$ in $\mathcal{G}_n$ converges
to $1$ as $n$ tends to infinity. Occasionally, we will say {\it
almost all} instead of almost every.

From the above theorem and employing Chebyshev inequality, it is
easy to see that
$$\mbox{Pr}\big[\big|X_n-E(X_n)\big|>
n^{3/4}\big]\leq \frac{Var X_n}{n^{3/2}}\rightarrow 0 \mbox{ as }
n\rightarrow \infty.$$ Thus, for almost all trees in $\mathcal
{T}_n$, $X_n$ equals $(\frac{2}{x_0b^2}\cdot w(i,j)+o(1))n$.
Consequently, the following result is relevant.

\begin{cor}\label{cor1} For almost all trees, the number of occurrences
of pattern $(i,j)$ is $(\frac{2}{x_0b^2}\cdot w(i,j)+o(1))n$.
\end{cor}

\section{An application}

In this section, we use the result of Corollary \ref{cor1} to
investigate the values of the {\it Randi\'{c} index} and {\it
general Randi\'{c} index}, and show that Conjecture \ref{rd} is true
for almost all trees.

Let $G=(V,E)$ be a graph with vertex set $V$ and edge set $E$. The
Randi\'{c} index is defined as
$$R(G)=\sum_{uv\in E} \frac{1}{\sqrt{d_ud_v}},$$
where $d_u$, $d_v$ are the degrees of the vertices $u, v\in V$.

We know that the number of occurrences of the pattern $(i,j)$ is the
number of edges with one end of degree $i$ and the other of degree
$j$ in the tree. Still, we assume $i\leq j$. Then, the number of
edges $(i,j)$ in almost all trees of $\mathcal{T}_n$ is
$(\frac{2}{x_0b^2}\cdot w(i,j)+o(1))n$. Moreover, every tree in
$T_n$ has $n-1$ edges. So, for any integer $K$, $\sum_{i\leq j\leq
K}(\frac{2}{x_0b^2}\cdot w(i,j))\leq 1$, it follows that
$\sum_{i\leq j}(\frac{2}{x_0b^2}\cdot w(i,j))$ is convergent.
Consequently, $\sum_{i\leq j}\frac{2}{x_0b^2\sqrt{i\cdot j}}\cdot
w(i,j)$ also converges to some constant $\lambda$. Although the
exact value of $\lambda$ can not be given, one can employ a computer
to get that $0.1<\lambda<1$. Then, for any $\varepsilon>0$, there
exists an integer $K_0$ such that for any $K\geq K_0$
$$\sum_{i\leq j, j\geq K}(\frac{2}{x_0b^2}\cdot w(i,j))<\varepsilon.$$
That is, for almost all trees, the number of edges
with one end of degree larger than $K$ is less than $\varepsilon n$.
Hence, the Randi\'{c} index enjoys
$$\big(\sum_{i\leq j\leq K}\frac{2}{x_0b^2\sqrt{i\cdot j}}w(i,j)+o(1)\big)\cdot n<R(T_n)
<\big(\sum_{i\leq j\leq K}\frac{2}{x_0b^2\sqrt{i\cdot
j}}w(i,j)+o(1)\big)\cdot n +\varepsilon n\mbox{ a.e.}$$ Immediately,
we obtain the following result.
\begin{thm} For any $\varepsilon>0$, the Randi\'{c} index of almost
all trees enjoys
\begin{eqnarray}\label{rg}
(\lambda-\varepsilon)n<R(T_n)<(\lambda+\varepsilon )n.
\end{eqnarray}
\end{thm}

Bollob{\'{a}s and Erd\"{o}s \cite{be} generalized the Randi\'{c}
index as
$$R_{\alpha}(G)=\sum_{uv\in E} (d_ud_v)^{\alpha},$$
which is called {\it general Randi\'{c} index}, where $\alpha$ is
a real number. Clearly, if $\alpha=-\frac{1}{2}$, then
$R_{-\frac{1}{2}}(G)=R(G)$. We refer the readers to a survey
\cite{lishi} for more details on this index. Here, we suppose
$\alpha< 0$. Following the sketch to obtain (\ref{rg}), we can
analogously get an estimate of $R_{\alpha}(T_n)$. Then, the
following result is relevant.

\begin{cor}
Suppose $\alpha<0$. For any $\varepsilon>0$ we have
$$(\lambda_{\alpha}-\varepsilon)n\leq R_{\alpha}(T_n)\leq(\lambda_{\alpha}+\varepsilon)n \mbox{    a.e.,}$$
where $\lambda_\alpha$ is some constant corresponding to every
$\alpha$.
\end{cor}

In what follows, we consider Conjecture \ref{rd}. Let $d(u,v)$ be
the distance between vertices $u, v\in V$. The {\it average
distance} is defined as the average value of the distances between
all pairs of vertices of a graph $G$, i.e.,
$$D(G)=\frac{\sum_{u,v\in V}d(u,v)}{{n\choose 2}}.$$

We will show that Conjecture \ref{rd} is true for almost all
trees. To this end, we first introduce the concept of {\it Wiener
index} for a graph $G$, which is defined as
$$W(G)=\sum_{u,v\in V}d(u,v).$$
Clearly, $W(G)={n\choose 2}D(G)$. $W(T_n)$ is a random variable on
$\mathcal {T}_n$, and Wagner \cite{w} established the following
result.

\begin{lem}
The Wiener index $W(T_n)$ enjoys
$$E(W(T_n))=(\omega+o(1))n^{5/2}$$ and
$$Var(W(T_n))=(\delta+o(1))n^{5},$$
where $\omega$ and $\delta$ are some constants.
\end{lem}

Employing Chebyshev inequality, we have
$$\mbox{Pr}[|W(T_n)-E(W(T_n))|\geq n^{11/4}]\leq\frac{Var(W(T_n))}{n^{11/2}}\rightarrow 0,
\mbox{ as } n\rightarrow\infty,$$ from the above lemma. Since
$E(W(T_n))=O(n^{5/2})$, therefore for almost all trees the Wiener
index $W(T_n)$ is $O(n^{11/4})$. Consequently, we can get that the
average distance satisfies
$$D(T_n)=O(n^{3/4})\mbox{ a.e.}$$
Combining with (\ref{rg}), the following
result is relevant.
\begin{thm}
For almost all trees in $\mathcal {T}_n$, $R(T_n)>D(T_n)$.
\end{thm}

\begin{rem} For the classic Erd\"{o}s--R\'{e}nyi model
$\mathcal{G}_{n,p}$ \cite{BB} of random graphs, which consists of
all graphs $G_{n,p}$ with vertex set $[n]=\{1,2,\ldots,n\}$ in which
the edges are chosen independently with probability $0<p<1$, we can
easily get the same result. In fact, suppose $p$ is a constant.
Recall that for almost all graphs the degree of a vertex is
$(p+o(1))n$ (see \cite{i}). Thus, for almost all graphs,
$$R(G_{n,p})=\frac{1}{2}\cdot\frac{1}{\sqrt{(p+o(1))^2n^2}}\cdot
(p+o(1))n\cdot n=(\frac{1}{2}+o(1))n.$$ Moreover, it is well known
that the diameter is not more than $2$ for almost all graphs.
Consequently, $D(G_{n,p})\leq 2 \mbox{  a.e.}.$ Hence,
$$R(G_{n,p})>D(G_{n,p})\mbox { a.e.}$$
\end{rem}

\noindent {\bf Acknowledgement:} The authors would like to thank
Dr. Gerard Kok for helps.

\end{document}